\documentclass[a4,11pt]{article}

\usepackage[english]{babel} 
\usepackage{amsmath,amssymb,amsthm,mathtools}
\usepackage{microtype}
\usepackage{enumitem}
\usepackage{hyperref}
\usepackage{tikz}
\usetikzlibrary{arrows.meta}
\usepackage{float}
\hypersetup{
colorlinks=true,
linkcolor=blue,
citecolor=blue,
urlcolor=blue
}

\newtheorem{theorem}{Theorem}[section]

\theoremstyle{definition}
\newtheorem{assumption}[theorem]{Assumption}

\theoremstyle{remark}

\newcommand{\LCM}{\mathcal S}

\newcommand{\pto}{\xrightarrow{P}}

\begin{document}

\title{Causality and order restricted inference in the Aalen additive hasard model  }
\author{Dragi Anevski\footnote{dragi@maths.lth.se} , Center for Mathematical Sciences, Lund University \\ and \\
ElBatoul Manel Merai\footnote{manel-elbatoul.merai@doc.umc.edu.dz}, Department of Mathematics, \\ Constantine 1 Brothers Mentouri University}\date{\today}
\date{}
\maketitle

\begin{abstract}
We study a graphical model (DAG) with time to event as a response variable. We use the Aalen additive hazard model for modelling the effect on the response variable of covariates, combined with a multiple linear regression for earlier layers in a DAG. This enables modelling covariates affecting the event process both directly and indirectly. The indirect component effect is represented by a product
$\theta_1(t)\beta_2(t)$, where $\theta_1(t)$ describes the effect of the independent covariate on the dependent covariate and $\beta_2(t)$ is the effect of the dependent covariate on the final time to response variable, in the additive hazard model. We estimate $\beta_2(t)$ by the order-restricted estimator recently introduced in \cite{AnevskiMerai2026}, thereby enabling disease models with early onset or late onset covariate effects. Combining our previously obtained results with the DAG allows for modelling early onset or late onset covariate effects together with direct and indirect covariate effects, which should be of particular interest for studying causality for complex diseases.
 
We present a nonparametric  plug-in estimator of the indirect effect. We derive the limit distribution of the  plug-in estimator of the indirect component, and show that it has  $n^{1/3}$ pointwise asymptotic rate and a limit distribution which is proportional to the Chernoff distribution.
\end{abstract}




\section{Introduction}

We consider a survival analysis model, in which we treat right censored data, combined with a DAG.  Thus we let $T\geq 0$ be a positive continuous random variable with an unknown distribution function $F$, where  $T$ models the time to an event of interest. We assume that we observe the minimum $t_i=\min(T_i,C_i)$ of the time $T_i$ and a censoring time $C_i$, together with an indicator for the time being exact, $\delta_i=1\{T_i\leq C_i\}$.  

Introduce the individual counting processes $N_i(t)=1\{t_i\leq t, \delta_i=1\}$ for which one has the stochastic differential equation
\begin{eqnarray}
   dN_i(t)&=& Y_i(t)h(t)dt+ dM_i(t), \label{eq:SDE-base}
\end{eqnarray}
for $i=1,\ldots,n$. Here $h(t)$ is the individual hasard function, which exists when the distribution of $T$ is absolutely continuous, and is then given by $h(t)=-{d}/{dt}\log(1-F(t))$, where $Y_i(t)=1\{t_i\geq t\}$ is the left-continuous indicator process for the individual being at risk at time $t-$. Also $M_i$ is the individual martingale, i.e. satisfying  $E(dM_i(t)| {\cal F}_t)=0$, and where the sigma algebras $\{{\cal F}_t, t\geq \}$ is a filtration storing the information available at times $t$.

The $\sigma$-algebra generated by the information depends on the amount and type of information about the observed times that is available and the amount of information about the covariates that is available. Starting with the model $(\ref{eq:SDE-base})$ satisfying $E(dM_i(t)| {\cal F}_t)=0$ one needs to establish that $E(dM_i(t)| {\cal G}_t)=0$ where ${\cal G}_t$ is the $\sigma$-algebra generated by the observed and available information at time $t$. The concept of noninformative and independent censoring as well as the innovation theorem are used to establish this link. We do not make this assumptions explicit here, since they are of no relevance to us, and refer to reader to a standard reference such as \cite{andersen.borgan.gill.keiding.1993}, cf.\;also  \cite{AnevskiMerai2026}. We will in the sequel assume that ${\cal F}_t$ is a filtration, containing the information available at time $t$, and that $E(dM_i(t)| {\cal F}_t)=0$ holds so that $E(dN_i(t)| {\cal F}_t)=Y_i(t)h(t)dt$.

A basic inference problem in survival analysis is assessing the effect of group indicators or continuous covariate measurements on the distribution  of the time to an event. It is then necessary to model the distribution function depending on covariates $x_1,\ldots,x_p$. This can be done using various models. One standard model which is particularly convenient in our framework is the Aalen additive hasard model
\begin{eqnarray}\label{eq:Aalen.model}
    h(t)&=&\beta_0(t)+\beta_1(t) x_1 +\ldots+\beta_p(t) x_p,
\end{eqnarray}
where $\beta_0,\ldots,\beta_p$ are unknown functions. Thus $\beta=(\beta_0,\ldots,\beta_p)$ is the unknown parameter vector of functions and $E^{p+1}$  is the parameter space, where $E=\{g:[0,\infty)\to {\mathbb R}, \int g<\infty \}$ is the set of integrable functions on $[0,\infty)$.  The Aalen model is immediately generalisable to time dependent covariates $x_k(t)$, which we will allow for in the sequel, and they may even be stochastic processes; the only necessary demand is that $X_k(t)$ is ${\cal F}_t$ measurable, for every $k$.

A principal difficulty in estimation in the Aalen model is that the coefficient functions are local functions and can not be estimated
directly by nonparametric procedures. Instead, one turns to estimating the cumulative
coefficient functions 
\begin{eqnarray}
B_k(t)=\int_0^t\beta_k(u)\,du. \label{eq:cumulative-functions}
\end{eqnarray}
The ordinary Aalen estimator, defined below in $(\ref{eq:aalen-B-estimator})$, provides an estimator $\widehat B_k$ of
$B_k$.

In \cite{AnevskiMerai2026} we proposed an order-restricted
estimator for the coefficient functions ${\beta}_k$ under a monotonicity assumption.
The estimator is obtained by projecting the local increments of the Aalen estimator $\widehat B_k$, defined in  $(\ref{eq:aalen-B-estimator})$, onto the cone of non-increasing sequences. Under some regularity conditions we showed in  \cite{AnevskiMerai2026}  that
the resulting pointwise estimator has a nonstandard $n^{1/3}$ asymptotic
rate and a Chernoff-type limit distribution.

One motivation for introduction order restricted inference of the $\beta_k$ functions is in order to model {\em early onset} covariate effects. In fact, we note that $\beta_k$ being decreasing means that the effect of (an increase in) the covariate value $x_k$ decreases with time or age.

A further motivation of the work in \cite{AnevskiMerai2026}, was as a first step in the use of a graphical model in a survival analysis setting, thereby allowing for a categorising of the covariate effects as {\em direct} versus {\em indirect} effect.  This was first suggested in \cite{fosen.borgan.weedon.aalen.2006}, see also \cite{aalen.borgan.gjessing.2008}. 

In this paper we incorporate our previously derived limit distribution results into such a graphical model setting and present estimators of momentanous direct and indirect effects, as well as limit distribution properties for the estimated effects. This has to our knowledge not been done before.

Thus the purpose of the present paper is to study a natural extension of a standard Aalen model, in which
one of the covariates acts through an intermediate, possibly time dependent, covariate. This leads to a decomposition of the effect  of a covariate into a direct component and an indirect component. For simplicity we will treat the simplest possible model, with two covariates. The generalisation to more complex models is straightfoward. 

More precisely, let $X_1(t)$ denote an (independent) covariate and let $X_2(t)$ denote another (dependent) covariate. We consider the structural representation
\begin{eqnarray}
dN_i(t) & = & Y_{i}(t)\{ \beta_0(t)+\beta_1(t)X_{1}(t) +  \beta_2(t)X_{2}(t) \}dt+dM_i(t),                                     
\label{eq:model}\\
X_{2}(t)
& = &
\theta_0(t)  + \ \theta_1(t)X_{1}(t)+\varepsilon(t),            \label{eq:covariate}
\end{eqnarray}
where $N_i$ is the event counting process, $Y_i$ is the at-risk process, $X_i(t)$ are ${\cal F}_t$-adapted processes and
$M_i$ is a martingale with respect to the filtration ${\cal F}_t$. 

The hazard model  is
\begin{eqnarray} \label{eq:hasard-representation}
h(t)& = & \beta_0(t)\ + \ \beta_1(t)X_{1}(t) \ + \ \beta_2(t)X_{2}(t).
\end{eqnarray}
Substitution of the covariate equation $(\ref{eq:covariate})$ into  $(\ref{eq:hasard-representation})$ gives
\begin{eqnarray*}
h(t)
& = &
\{\beta_0(t)+\beta_2(t)\theta_0(t)\} + 
\{\beta_1(t)+\theta_1(t)\beta_2(t)\}X_{1}(t)
+ 
\beta_2(t)\varepsilon(t).
\end{eqnarray*}

\begin{figure}[H]
\centering
\begin{tikzpicture}[>=Stealth, node distance=2.5cm,
every node/.style={font=\small}]
\node (x1) at (0,0) {$X_1(t)$};
\node (x2) at (3,0) {$X_2(t)$};
\node (dn) at (1.5,-2) {$dN(t)$};

\draw[->] (x1) -- (x2)
node[midway, above] {$\theta_1(t)$};

\draw[->] (x1) -- (dn)
node[midway, left, xshift=-2pt] {$\beta_1(t)$};

\draw[->] (x2) -- (dn)
node[midway, right, xshift=2pt] {$\beta_2(t)$};
\end{tikzpicture}
\caption{Dynamic path diagram illustrating the direct and indirect pathways from $X_1(t)$ to the event process $N(t)$.}
\label{fig:diagram}
\end{figure}
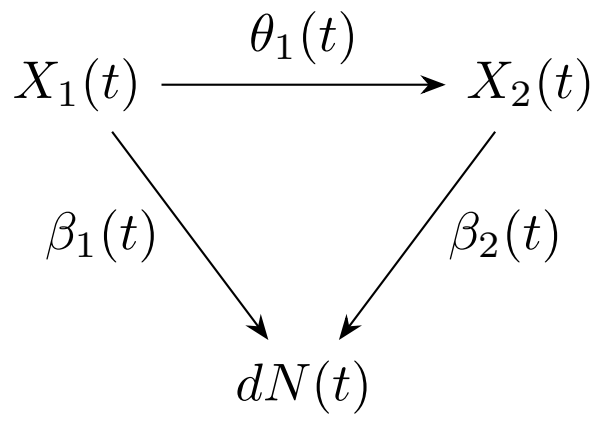

Thus, at the level of this structural decomposition, the coefficient
multiplying $X_{1}(t)$ is
\begin{eqnarray}
\underbrace{\beta_1(t)}_{\text{direct component}}
 &+&  
\underbrace{\theta_1(t)\beta_2(t)}_{\text{indirect component}}.
\label{eq:decomp}
\end{eqnarray}
The two pathways represented in \eqref{eq:decomp} are illustrated in Figure~\ref{fig:diagram}.

We thus define the instantaneous {\em indirect effect} of covariate $X_1(t)$ as 
\begin{eqnarray*}
I(t) & = & \theta_1(t)\beta_2(t).
\end{eqnarray*}
Our main question is the point wise asymptotic behaviour of the plug-in estimator 
\begin{eqnarray*}
\widehat I_n(t) & = & \widehat\theta_1(t)\widetilde\beta_2(t)
\end{eqnarray*}
of the instantaneous indirect effect $I(t)$, where $\widetilde\beta_2$ is the order restricted estimator defined in \cite{AnevskiMerai2026}, and where $\widehat\theta_1$ is a least squares estimator, to be defined in the sequel.

We will show below that the two factors $\widehat\theta_1$ and $\widetilde\beta_2$ have different asymptotic
rates. In fact, for each fixed $t$, the estimator $\widehat\theta_1(t)$ is obtained
by ordinary least squares regression and satisfies a standard
$n^{1/2}$ distributional limit, whereas $\widetilde\beta_2(t)$ has a nonstandard
$n^{1/3}$ limit distribution behaviour that is common in order restricted inference. Since
$n^{-1/2}=o(n^{-1/3})$, estimation of $\theta_1(t)$ does not contribute to the first order asymptotic distribution of the estimated  indirect effect $\widehat I_n(t)$.

The main result therefore transfers the Chernoff-type limit of
$\widetilde\beta_2(t)$ directly to the indirect component, with the limit
multiplied by the fixed value $\theta_1(t)$. The proof is straightforward once
the pointwise result for the order restricted estimator is available, but
the formulation is useful because it identifies precisely which part of
the estimation procedure determines the asymptotic behaviour.

The paper is organised as follows. Section~2 introduces the model and the
estimators. Section~3 derives the limit distribution for the new estimator. Section~4 discusses the interpretation and
limitations of the result.

\section{Model and estimation}

Assume we have data $(t_i,\delta_i,x_{1i}(t),x_{2i}(t)$ for $i=1,\ldots,n$ for $n$ individuals, which are seen as observations in the counting process Aalen model $(\ref{eq:model})$, which satisfies the structural equations $(\ref{eq:covariate})$. Note that we allow the observed covariates $x_{1i},x_{2i}$ to be time-dependent and arbitrary. The only requirement needed is that they are observations of ${\cal F}_t$-adapted processes, so that the martingale structure is maintained under stochastic integrations. We assume we have observations on a fixed compact time interval $[0,\tau]$.


Now for each time point $t\in [0,\tau]$, let $J(t)=\{i:Y_i(t)=1\}$ be the {\em risk set}, i.e.\,the index set of those individuals that have not yet experienced the event.
Since we do not have left-truncated data, it is clear that $J(t)$ is a set of indices that is decreasing with time, so that $J(t)\supseteq J(s)$ if $t<s$. Furthermore  the cardinality of $J(t)$ decreases with one at each observation time $t_i$, so that $|J(t_i+)|=|J(t_i)|-1$, assuming that we have no ties, which we hereby assume.

Therefore, for each fixed observed time point $t_i$, for $i=1,\ldots,n$, the structural equations $(\ref{eq:covariate})$ give rise to a set of linear equations 
\begin{eqnarray}
X_{2j}(t_i) & = &\theta_0(t_i)  +  \theta_1(t_i)X_{1j}(t_i)+\varepsilon_j(t_i),
\label{eq:med}
\end{eqnarray}
for $j\in J(t_i)$, and we note that there are $n_{i}=|J(t_i)|$ many such equations. 

Therefore, for each fixed $t_i$, we can estimate the parameters $\theta_0(t_i),\theta_1(t_i)$ by ordinary least squares (ols) regression. We denote the ols estimates  by $(\hat{\theta}_0(t_i),\hat{\theta}_1(t_i))$. Since we can solve the system of equations $(\ref{eq:med})$ for every time point $t_i$, we will obtain nonparametric estimates of the full parameter function $\theta_0(t)$ and $\theta_1(t)$, and these are
\begin{eqnarray*}
     \hat{\theta}_0&:=&(\hat{\theta}_0(t_{(1)}),\ldots,\hat{\theta}_0(t_{(n)})),\\
     \hat{\theta}_1&:=&(\hat{\theta}_1(t_{(1)}),\ldots,\hat{\theta}_1(t_{(n)})).
\end{eqnarray*}
In fact $ \hat{\theta}_0$ and $\hat{\theta}_1$ are the ols estimate of the vectors $\theta_0:=(\theta_0(t_{(1)},\ldots, \theta_0(t_{(n)})$ and $ \theta_1:=(\theta_1(t_{(1)},\ldots, \theta_1(t_{(n)})$. 

Note that even if the covariates would be time-constant, the ols estimates  $(\hat{\theta}_0(t_i),\hat{\theta}_1(t_i))$ will typically be different for different $t_i$, since the risk set $J(t_i)$ changes with $t_i$, and thus the set of linear equations $(\ref{eq:med})$ will change as well (the number of equations will change). 

Note also that if the covariates would be time-dependent then we could formulate system of equation $(\ref{eq:med})$ for every fixed $t \in (0,\tau)$, with $t_i$ replaced by $t$ and the index set $J(t_i)$ replaced by $J(t)$. We would then get ols estimates of $\hat{\theta}_0(t)$ and $\hat{\theta}_1(t)$ of $\theta_0(t)$ and $\theta_1(t)$. Since we can do this for any fixed $t$ we actually get nonparametric estimates of the parameter functions, under very mild regularity conditions on the covariates, e.g.\,that they are piecewise continuous. 

Then a multivariate central limit theorem will give a limit distribution result for the $2n$-length vector $( \hat{\theta}_0,\hat{\theta}_1)$ as
\begin{eqnarray*}
n^{1/2}( (\hat{\theta}_0,\hat{\theta}_1) - (\theta_0,\theta_1))
&\xrightarrow{d}& U,
\end{eqnarray*}
on ${\mathbb R}^{2n}$ where $U$ is a mean zero $2n$-dimensional  Gaussian random vector, with a certain covariance structure. It is also straightforward to show that we have a process result
 \begin{eqnarray*}
n^{1/2}( (\hat{\theta}_0(\cdot),\hat{\theta}_1(\cdot)) - (\theta_0(\cdot),\theta_1(\cdot)))
&\xrightarrow{d}& U(\cdot),
\end{eqnarray*}
on $D[0,\tau]^2$, where $U(\cdot)$ is a mean zero Gaussian process in $D[0,\tau]^2$, with a certain covariance function, under mild regularity conditions on the covariates.

For the purpose of the results presented in this paper, only the limit distribution of $\hat{\theta}_1(t_0)$ at a fixed point interior $t_0\in (0,\tau)$ is of interest. We state that limit distribution as an assumption, to allow for general statements.
\begin{assumption} \label{eq:thet}
For the fixed interior point $t_{0}\in (0,\tau)$, the ols estimator
$\widehat{\theta}_1(t_{0})$ satisfies
\begin{eqnarray}
n^{1/2}\left(\widehat\theta_1(t_{0})-\theta_1(t_{0})\right)
\xrightarrow{d} G(t_{0}),
\end{eqnarray}
where $G(t_0)$ is a mean zero random variable.
\end{assumption}
Note that this assumption is satisfied for the case that we treat above of independent data, with $G(t_0)$ a mean zero Gaussian random variable, under mild regularity conditions.

Now, recall the definition of the cumulative functions $B_k$, for $k=0,1,2$,  in $(\ref{eq:cumulative-functions})$. Then a standard result, cf. e.g.\, \cite{andersen.borgan.gill.keiding.1993}, is that the least squares estimator of $B=(B_0,B_1,B_2)^T$ is
\begin{eqnarray}
\widehat B(t)
& = &
\int_0^t In(u)Y^{-}(u)\,dN(u), \label{eq:aalen-B-estimator}
\end{eqnarray}
where $In(u)$ is a indicator for a certain design matrix in a matrix formulation of $(\ref{eq:SDE-base})$ and $(\ref{eq:Aalen.model})$ has full rank, and $Y^{-}(u)$
denotes a generalized inverse of the corresponding design matrix; see
\cite{Aalen1980} and \cite{andersen.borgan.gill.keiding.1993} for more details and \cite{AnevskiMerai2026} for a review.

Recently, in \cite{AnevskiMerai2026}, this estimator was used as an input into an order-restricted scheme, to obtain a monotone estimator of the $\beta_k$ functions themselves. Thus if we assume that $\beta_2$ is non-increasing, and $\widehat B_2$ denote the
second coordinate of $\widehat B$, then \cite{AnevskiMerai2026} defined an order-restricted estimator
$\widetilde\beta_2$, namely the isotonic regression of the local
increment estimates of $\widehat B_2$ onto the cone of non-increasing
sequences. Furthermore, in \cite{AnevskiMerai2026} it was shown that a characterisation of the solution is
\begin{eqnarray}
\widetilde\beta_2(t) & = &
\frac{d}{dt}\LCM{\widehat B_2}(t),
\label{eq:lcm}
\end{eqnarray}
where the derivative is understood to be the appropriate one sided derivative. Furthermore, algorithms for calculating the isotonic regression are well known, and e.g.\,the minimum lower set algorithm works well, cf.\,e.g.\,\cite{Robertson1988}.

Now, since we have obtained estimators of the factors $\beta_2$ and $\theta_1$ of the indirect effect $I(t)$, we can define the plug-in estimator of the indirect effect  by
\begin{eqnarray}
\widehat I_n(t) & = &
\widehat\theta_1(t)\widetilde\beta_2(t). \label{eq:pluin-estimator-def}
\end{eqnarray}

\section{Limit distributions for the indirect effect estimator}

We first recall and state the assumptions for the survival analysis setting that enable the limit distribution results  of \cite{AnevskiMerai2026} for the order restricted estimator  $\widetilde\beta_2(t)$. In particular, these conditions ensure local
regularity of the design and risk processes, nonsingularity of the limiting
design matrix, and the smoothness and strict local monotonicity conditions
on the coefficient function at the point under consideration. We emphasise that they are equivalent to the standard assumptions for deriving limit distributions of the Aalen estimator $(\ref{eq:aalen-B-estimator})$, cf.\,\cite{andersen.borgan.gill.keiding.1993}, and thus we do not need more than what is regularly assumed in a standard survival analysis setting.

\begin{assumption}\label{eq:ct}
Assume that Assumptions 1,2 and 3 in \cite{AnevskiMerai2026} hold and that  $\beta_2$ is differentiable with $\beta_2'<0$ in a neighbourhood around the fixed interior point $t_0\in (0,\tau)$. 
\end{assumption}
Under Assumption \ref{eq:ct}, Theorem 3 in \cite{AnevskiMerai2026}
gives the limit distribution
\begin{eqnarray}
n^{1/3}c(t_0)
\left(\widetilde\beta_2(t_0)-\beta_2(t_0)\right)
\xrightarrow{d} \LCM\left(-s^2+w(s)\right)'(0)
\label{T}
\end{eqnarray}
where  $w$ is a standard two-sided Brownian motion, and 
 where 
  \begin{eqnarray*}
         c(t_0)&=& 2^{-1/3}|\beta_2'(t_0)|^{1/3} (\sigma_2^2)^{-4/3},
  \end{eqnarray*}
and, $\sigma_2^2$ is specified in \cite{AnevskiMerai2026}.

Next we state the main result on the limit distribution for the plug estimator $\widehat I_n(t)$, defined in $(\ref{eq:pluin-estimator-def})$, of the indirect effect $I(t)$.

\begin{theorem}
Suppose that Assumptions \ref{eq:thet} and \ref{eq:ct} hold. Then
\begin{eqnarray*}
n^{1/3}c(t_{0})
\left(\widehat I_n(t_{0})-I(t_{0})\right)
\xrightarrow{d}
\theta_1(t_{0})\LCM\left(-s^2+w(s)\right)'(0).
\end{eqnarray*}
\end{theorem}
We note that the limit distribution is proportional to the Chernoff random variable $\LCM\left(-s^2+w(s)\right)'(0)$, common in order restricted inference. 
\begin{proof}
We have that 
\begin{eqnarray*}
n^{1/3}
\left(\widehat I_n(t_{0})-I(t_{0})\right)
&=
\widehat\theta_1(t_{0})
n^{1/3}
\left(\widetilde\beta_2(t_{0})-\beta_2(t_{0})\right)\\ \nonumber
& + 
\beta_2(t_{0})n^{1/3}
\left(\widehat\theta_1(t_{0})-\theta_1(t_{0})\right).
\end{eqnarray*}
We first see that the second term is
\begin{eqnarray*}
n^{1/3}\left(\widehat\theta_1(t_{0})-\theta_1(t_{0})\right)
& = &
n^{-1/6}
\sqrt n\left(\widehat\theta_1(t_{0})-\theta_1(t_{0})\right).
\end{eqnarray*}
and by Assumption \ref{eq:thet}, since $n^{-1/6}\to0$ and the second factor is bounded in probability,
the product converges to zero in probability, and therefore the second term is  $o_P(1)$.

Next we treat the first term: Since $\widehat\theta_1(t_{0})-\theta_1(t_{0})\pto0$, and since by \eqref{T} $n^{1/3}
\left(\widetilde\beta_2(t_{0})-\beta_2(t_{0})\right)=  O_P(1)$, we have that
\begin{eqnarray*}
\left(\widehat\theta_1(t_{0})-\theta_1(t_{0})\right)
n^{1/3}
\left(\widetilde\beta_2(t_{0})-\beta_2(t_{0})\right)
& = & o_P(1),
\end{eqnarray*}
Therefore
\begin{eqnarray*}
\widehat\theta_1(t_{0}) n^{1/3}\big(\widetilde\beta_2(t_{0})-\beta_2(t_{0})\big) & = &
\theta_1(t_{0}) n^{1/3}\big(\widetilde\beta_2(t_{0})-\beta_2(t_{0})\big)
\\
&& +  \big(\widehat\theta_1(t_{0})-\theta_1(t_{0})\big) n^{1/3}\big(\widetilde\beta_2(t_{0})-\beta_2(t_{0})\big)\\
& = &
\theta_1(t_{0})
n^{1/3}
\left(\widetilde\beta_2(t_{0})-\beta_2(t_{0})\right)
\ + \ o_P(1).
\end{eqnarray*}

Thus

\begin{eqnarray*}
n^{1/3}
\left(\widehat I(t_{0})-I(t_{0})\right)
& = &
 \theta_1(t_{0})
n^{1/3}
\left(\widetilde\beta_2(t_{0})-\beta_2(t_{0})\right)
\ + \ o_P(1).
\end{eqnarray*}

Multiplying by the fixed constant $c(t_{0})$ and applying
the limit distribution result \eqref{T} gives

\begin{eqnarray*}
n^{1/3}c(t_{0})\left(\widehat I_n(t_{0})-I(t_{0})\right)
\xrightarrow{d}
\theta_1(t_{0})\LCM\left(-s^2+w(s)\right)'(0).
\end{eqnarray*}

\end{proof}

\section{Conclusion and discussion}

We have established the pointwise asymptotic distribution of a plug-in estimator for the indirect component in an DAG with the final stage modeled by an additive hazard model. When the coefficient $\theta_1$ modelling the independent to dependent covariate effect is estimated at the regular $ n^{1/2}$ rate and the covariate coefficient in the final Aalen model is estimated by the newly developed order restricted Aalen procedure in \cite{AnevskiMerai2026}, we show that the latter determines the asymptotic behaviour of the plug-in estimator of the indirect effect. In particular the resulting plug-in estimator has the same $n^{1/3}$ pointwise scale as the order-restricted coefficient, with a Chernoff-type limit multiplied by the exposure to covariate coefficient.

The result provides a theoretical basis for inference on an indirect component constructed from an order-restricted additive hazard estimator.

The suggested model and estimator is novel, and allows for studying both {\em early onset} covariate effects and {\em indirect} covariate effects, and should be highly interesting for modelling complex diseases with time-to-event as response variable.


\end{document}